\documentclass[secthm,seceqn,amsthm,ussrhead,12pt]{amsart}
\usepackage{amsmath,latexsym}
\usepackage[english]{babel}
\usepackage[psamsfonts]{amssymb}
\usepackage{times}
\usepackage{cite}

\usepackage[mathcal]{euscript}

\begin{document}
	\sloppy
	
\begin{center}
{\huge Generalized derivations of (color) $n$-ary  algebras.}

\medskip

\medskip
\textbf{Ivan Kaygorodov$^{a,b,c}$, Yury Popov$^{c,d}$}\footnote{The authors were supported
by RFBR 15-51-04099.
The first author was supported by FAPESP 14/24519-8.}
\medskip

$^a$ Universidade Federal do ABC, CMCC, Santo Andr\'{e}, Brasil,\\
$^b$ Universidade de Sao Paulo, Sao Paulo,  Brasil,\\
$^c$ Sobolev Institute of Mathematics, Novosibirsk, Russia,\\
$^d$ Novosibirsk State University, Russia.

\end{center}

       \vspace{0.3cm}

{\bf Abstract.}
We generalize the results of Leger and Luks about generalized derivations of Lie algebras
to the case of color $n$-ary $\Omega$-algebras.
Particularly, we prove some properties of generalized derivations of color $n$-ary algebras;
prove  that a quasiderivation algebra of a color  $n$-ary $\Omega$-algebra can be embedded into
the derivation algebra of a larger color  $n$-ary $\Omega$-algebra,
and describe (anti)commutative $n$-ary algebras satisfying the condition $QDer = End.$

       \vspace{0.3cm}

\noindent {\bf \S 0\quad Introduction}

       \vspace{0.3cm}

It is well known that the algebras of derivations and generalized derivations are very important in the study of Lie algebras and its generalizations.
There are many generalizations of derivations
(for example, Leibniz derivations \cite{KP} and Jordan derivations \cite{her}).
The notion of $\delta$-derivations appeared in the paper of Filippov \cite{fil1},
he studied $\delta$-derivations of prime Lie and Malcev algebras \cite{fil2,fil3}.
After that, $\delta$-derivations of
Jordan and Lie superalgebras were studied in \cite{kay1, kay_lie, kay_lie2,kay12mz, zus10} and many other works.
The notion of generalized derivation is a generalization of $\delta$-derivation.
The most important and systematic research on the generalized derivations algebras of a Lie algebra and their subalgebras was
due to Leger and Luks \cite{LL00}.
In their article, they studied properties of generalized derivation algebras and their subalgebras, for example, the quasiderivation algebras. They have determined the structure of algebras of quasiderivations and generalized derivations and proved that the     quasiderivation algebra of a Lie algebra can be embedded into the derivation algebra of a larger Lie algebra.
Their results were generalized by many authors. For example, Zhang and Zhang \cite{ZZ10} generalized the above results to the case of Lie superalgebras;
Chen, Ma, Ni and Zhou considered the generalized derivations of color Lie algebras, Hom-Lie superalgebras and Lie triple systems \cite{cml13,cml14,cml15}.
Generalized derivations of simple algebras and superalgebras
were investigated in \cite{KN03,lesha12,lesha14,GP03}.
Perez-Izquierdo and Jimenez-Gestal used the generalized derivations to study non-associative algebras \cite{GP08,GP09}.
Derivations and generalized derivations of $n$-ary algebras were considered in \cite{poj06,ck10,kay12izv,kay11aa,kay14sp,kay14mz,Williams} and other.
For example, Williams proved that, unlike the case of binary algebras, for any $n \geq 3$ there exist a non-nilpotent $n$-Lie algebra with invertible derivation  \cite{Williams},
Kaygorodov described $(n+1)$-ary derivations of simple $n$-ary Malcev algebras \cite{kay14sp} and  generalized derivations algebra of semisimple Filippov algebras over an algebraically closed field of characteristic zero \cite{kay14mz}.

The main purpose of our work is to generalize the results of Leger and Luks to the case of color $n$-ary algebras.
Particularly, we prove some properties of generalized derivations of color $n$-ary algebras;
prove  that the quasiderivation algebra of a color  $n$-ary $\Omega$-algebra can be embedded into
the derivation algebra of a large color  $n$-ary $\Omega$-algebra.
We describe all nonabelian $n$-ary (anti)commutative algebras with the condition $QDer = End.$
Particularly, we prove that if a Filippov algebra satisfies this property then it is either the solvable $n$-dimensional algebra or the simple $(n+1)$-dimensional algebra.

       \vspace{0.3cm}


\noindent {\bf \S 1\quad  Preliminaries}

        \vspace{0.3cm}

Let $F$ be a field and $G$ be an abelian group.
A map $\epsilon: G \times G \rightarrow F^*$ is called a bicharacter on $G$ if the
following relations hold for all $f,g,h \in G:$

(1) $\epsilon(f,g+h)=\epsilon(f,g)\epsilon (f,h),$

(2) $\epsilon(g+h,f)=\epsilon(g,f)\epsilon(h,f),$

(3) $\epsilon(g,h)\epsilon(h,g)=1.$
\vspace{0.3cm}

 \noindent{\bf Definition 1.1.} {
A color  $n$-ary algebra $T$ is an $n$-ary $G$-graded vector space $T=\bigoplus\limits_{g\in G} T_g$ with a graded $n$-linear map
$[\cdot , \ldots, \cdot ]: T \times \ldots \times T \rightarrow T$ satisfying

$$[T_{\theta_1},  \ldots, T_{\theta_n}] \subseteq T_{\theta_1+ \ldots +\theta_n}, \forall \theta_i \in G.$$
The main examples of color $n$-ary algebras
are color Lie algebras \cite{cml13,mikh85},
color Leibniz algebras \cite{dzhuma},
Filippov ($n$-Lie) superalgebras \cite{poj08patricia,ck10,poj03,poj09,poj08}
and $3$-Lie colour algebras \cite{ZT15}.

Let $T=\bigoplus\limits_{g\in G} T_g$ be a color algebra.
An element $x$ is called a homogeneous element of degree $t \in G$ if $x \in T_t$. We denote this by $hg(x)=t$.
A linear mapping $D$ is homogeneous of degree $t$ if $\forall g \in G: D(T_g) \subseteq T_{g+t}$
We denote this by $hg(D)=t.$ From now on, unless stated otherwise we assume that all elements and mappings are homogeneous. For two homogeneous elements $a$ and $b$, we set $\epsilon(a,b):=\epsilon(hg(a),hg(b)).$

We use the following notations:
\begin{eqnarray}\label{colorcom}[D_1,D_2] = D_1D_2 - \epsilon(D_1,D_2)D_2D_1,
\end{eqnarray}
\begin{eqnarray*}X_i = hg(x_1)+\ldots +hg(x_i).
\end{eqnarray*}

By ${\rm End}(T)$ we denote the set of all linear maps of
$T$. Obviously, ${\rm End}(T)$ endowed with the color bracket (\ref{colorcom}) is a color Lie algebra over $\mathbb{F}$.
\vspace{0.3cm}

\noindent{\bf Definition 1.2.}  {\it Let $(T,[\cdot, \ldots ,\cdot])$ be a color $n$-ary algebra.
A linear map
$D:T\rightarrow T$ is said to be a derivation of $T$ if it satisfies
$$ \sum \epsilon(D, X_{i-1}) [x_1, \ldots, D(x_i), \ldots, x_n]=D([x_1, \ldots, x_n]),$$
for all $x_1, \ldots, x_n \in T.$}

We denote the set of derivations of $T$ by ${\rm Der}(T)$. It is easy to see that ${\rm
Der}(T)$ endowed with the usual color commutator is a subalgebra of ${\rm
End}(T)$. The algebra ${\rm Der}(T)$ is called the derivation algebra of $T$.\vspace{0.3cm}

\noindent{\bf Definition 1.3} {\it A linear mapping $D\in {\rm End}(T)$ is called a generalized
 derivation of degree $\theta$ of $T$ if there exist linear mappings
$D',D'',\ldots, D^{(n-1)},D^{(n)} \in {\rm End}(T)$ of degree $\theta$ such that
$$\sum \epsilon(D, X_{i-1}) [x_1, \ldots, D^{(i-1)}(x_i), \ldots, x_n]=D^{(n)}([x_1, \ldots, x_n])$$
for all $x_1, \ldots, x_n \in T.$ An $(n+1)$-tuple $(D, D', \ldots, D^{(i-1)},  \ldots, D^{(n-1)},D^{(n)})$ is called an $(n+1)$-ary derivation.
} \vspace{0.3cm}

\noindent{\bf Definition 1.4.} {\it A linear mapping $D\in {\rm End}(T)$ is said to be a quasiderivation of degree $\theta$ if there exists a $D'\in
{\rm End}(T)$ of degree $\theta$ such that
$$\sum \epsilon(D, X_{i-1}) [x_1, \ldots, D(x_i), \ldots, x_n] =D^{'}([x_1, \ldots, x_n])$$
for all $x_1, \ldots, x_n \in T.$}

The sets of generalized derivations and quasiderivations will be denoted by ${\rm GDer}(T)$ and ${\rm QDer}(T)$
respectively.\vspace{0.3cm}

\noindent{\bf Definition 1.5.}
{\it The set ${\rm C(T)}$ consisting of linear mappings $D$ with the property
$$\epsilon(D, X_{i-1}) [x_1, \ldots, D(x_i),\ldots, x_n ]= D([x_1, \ldots, x_n]) \text{ for all } x_1, \ldots, x_n \in T$$
is called the centroid of $T$.}\vspace{0.3cm}

 \noindent{\bf Definition 1.6.}
{\it The set ${\rm QC(T)}$ consisting of linear mappings $D$ with the property
$$[D(x_1), \ldots, x_n]=\epsilon(D, X_{i-1})[x_1, \ldots, D(x_i), \ldots, x_n] \text{ for all } x_1, \ldots, x_n \in T$$
is called the quasicentroid of $T$.}\vspace{0.3cm}

\noindent{\bf Definition 1.7.} {\it The set ${\rm ZDer(T)}$ consisting of linear mappings $D$ with the property
$$[x_1, \ldots, D(x_i), \ldots, x_n]=D([x_1, \ldots, x_n])=0 \text{ for all } x_1, \ldots, x_n \in T, i = 1, \dots, n$$
is called the set of central derivations of $T$.}

It is easy to verify that $${\rm ZDer}(T)\subseteq {\rm
Der}(T)\subseteq {\rm QDer}(T)\subseteq {\rm GDer}(T)\subseteq {\rm
End}(T).$$



       \vspace{0.3cm}

     \noindent{\bf \S 2\quad Generalized derivation algebras and their color subalgebras}

       \vspace{0.3cm}

  First, we  give some basic properties
of the center derivation algebra, quasiderivation algebra and the
generalized derivation algebra of a color $n$-ary algebra.

     \noindent{\bf  Lemma 2.1.} {\it Let
$T$ be a color $n$-ary algebra. Then the following statements hold:

$(1)$\quad ${\rm GDer}(T),{\rm QDer}(T)$ and ${\rm C}(T)$ are color
subalgebras of  ${\rm End}(T)$;

$(2)$\quad ${\rm ZDer}(T)$ is a color ideal of ${\rm Der}(T)$.}

\rm {\it Proof.}\quad  $(1)$ Let $D_{1},D_{2}\in{\rm
GDer}(T)$. For arbitrary $x_1, \ldots, x_n \in T$ we have

$$[D_{1}D_{2}(x_1),\ldots, x_n] =
D_1^{(n)}([D_2(x_1), \ldots , x_n]) - \sum\limits_{i=2}^n \epsilon(D_1, D_2+X_i) [D_2(x_1), \ldots, D_{1}^{(i-1)}(x_i), \ldots, x_n]=$$
$$D_1^{(n)}D_2^{(n)}([x_1, \ldots, x_n]) - \sum_{j=2}^n  \epsilon(D_2, X_j)D_1^{(n)}([ x_1, \ldots, D_2^{(j-1)}(x_j), \ldots, x_n])-$$
$$\sum\limits_{i=2}^n \epsilon(D_1, D_2+X_i) D_2^{(n)}([x_1, \ldots, D_1^{(i-1)}(x_i), \ldots, x_n]) +$$
$$\sum\limits_{i=2,j=2, j< i}^n \epsilon(D_2,X_j)\epsilon(D_1, D_2+X_i) [x_1, \ldots, D_2^{(j-1)}(x_j), \ldots, D_1^{(i-1)}(x_i), \ldots, x_n]+$$

$$\sum\limits_{i=2,j=2, i<j}^n \epsilon(D_1,D_2+X_i)\epsilon(D_2, D_1+X_j) [x_1, \ldots, D_1^{(i-1)}(x_i), \ldots, D_2^{(j-1)}(x_j), \ldots, x_n]+$$
$$ \sum\limits_{i=2}^n \epsilon(D_1,X_i)\epsilon(D_2, X_i) [x_1, \ldots, D^{(i-1)}_1D_2^{(i-1)}(x_i), \ldots, x_n].$$

\noindent Thus for arbitrary $x_1,\ldots, x_n \in T$ we have
$$[[D_{1},D_{2}](x_1), \ldots, x_n]=[(D_{1}D_{2}- \epsilon(D_1,D_2)D_2D_1)(x_1), \ldots, x_n]=$$
$$[D^{(n)}_{1},D^{(n)}_{2}]([x_1, \ldots, x_n])-
\sum\limits_{i=2}^n \epsilon(D_1+D_2, X_i)[x_1, \ldots, [D^{(i-1)}_{1},D^{(i-1)}_{2}](x_i), \ldots, x_n].$$

From the definition of generalized derivation one gets
$[D_{1},D_{2}]\in {\rm GDer}(T),$ so ${\rm GDer}(T)$ is a color subalgebra of ${\rm End}(T)$.

Similarly, ${\rm QDer}(T)$ is a color subalgebra of ${\rm End}(T)$.

Let $D_{1},D_{2}\in{\rm C}(T).$ Therefore, for arbitrary $x_1, \ldots, x_n \in T$ we have
$$[[D_{1},D_{2}](x_1), \ldots, x_n]=[D_{1}D_{2}(x_1), \ldots, x_n]-\epsilon(D_1,D_2)[D_{2}D_{1}(x_1), \ldots, x_n]=$$
$$D_{1}([D_{2}(x_1), \ldots, x_n])-\epsilon(D_1,D_2)D_{2}([D_{1}(x_1), \ldots, x_n])=$$
$$D_{1}D_{2}([x_1, \ldots, x_n])-\epsilon(D_1,D_2)D_{2}D_{1}([x_1, \ldots, x_n])=$$
$$[D_{1},D_{2}]([x_1, \ldots, x_n]).$$

Similarly,
$$ \epsilon(D_1+D_2, X_i)[x_1, \ldots, [D_{1},D_{2}](x_i), \ldots, x_n]=[D_{1},D_{2}]([x_1, \ldots, x_n]).$$
Then $[D_{1},D_{2}]\in {\rm C}(T)$ and ${\rm C}(T)$ is a color subalgebra of
${\rm End}(T)$.

$(2)$ Let $D_{1}\in{\rm ZDer}(T),D_{2}\in{\rm Der}(T).$ For
arbitrary $x_1, \ldots, x_n \in T$ we have
$$ \epsilon(D_1+D_2, X_i) [x_1, \ldots, [D_{1},D_{2}](x_i), \ldots, x_n]=$$
$$\epsilon(D_1+D_2, X_i) [x_1, \ldots, (D_{1}D_{2}-\epsilon(D_1,D_2)D_{2}D_{1})(x_i),\ldots, x_n]=$$

$$\epsilon(D_1+D_2, X_i) [x_1, \ldots, D_{1}(D_{2}(x_i)),\ldots, x_n] -
\epsilon(D_1,D_2)\epsilon(D_1, X_i) D_2( [x_1, \ldots, D_{1}(x_i),\ldots, x_n]) + $$
$$\sum\limits_{j=1}^{i-1} \epsilon(D_1,D_2)\epsilon(D_2,X_j)  \epsilon(D_1, X_i) [x_1, \ldots, D_2(x_j), \ldots, D_1(x_i) ,\ldots, x_n]+$$
$$\sum\limits_{j=i+1}^n \epsilon(D_1,X_i)\epsilon(D_2,X_j) [x_1, \ldots, D_1(x_i), \ldots, D_2(x_j) ,\ldots, x_n]=0.$$

It is easy to see that $[D_1,D_2]$ is a derivation and
$$[D_1,D_2] ([x_1, \ldots, x_n])=0.$$
Then $[D_{1},D_{2}]\in {\rm ZDer}(T)$ and ${\rm ZDer}(T)$ is an
ideal of ${\rm Der}(T)$.\hfill$\Box$\vspace{0.3cm}


\noindent{\bf Lemma 2.2.} {\it Let $T$  be a color $n$-ary algebra. Then

    $(1)$ \quad $[{\rm Der}(T),{\rm C}(T)]\subseteq {\rm C}(T);$

    $(2)$ \quad $[{\rm QDer}(T),{\rm QC}(T)]\subseteq {\rm QC}(T);$

    $(3)$ \quad ${\rm C}(T)\cdot{\rm Der}(T)\subseteq {\rm Der}(T);$

    $(4)$ \quad ${\rm C}(T)\subseteq {\rm QDer}(T);$

    $(5)$ \quad $[{\rm QC}(T),{\rm QC}(T)]\subseteq {\rm QDer}(T);$

    $(6)$ \quad ${\rm QDer}(T)+{\rm QC}(T)\subseteq {\rm GDer}(T).$}

 \rm {\it Proof.}\quad $(1)$ Let
$D_{1}\in{\rm Der}(T),D_{2}\in{\rm C}(T).$ For arbitrary $x_1, \ldots, x_n \in T$ we
have
$$ [D_{1}D_{2}(x_1), \ldots, x_n]=D_{1}([D_{2}(x_1), \ldots, x_n])-
\sum\limits_{i=2}^n \epsilon(D_1, D_2+X_i)[D_{2}(x_1),\ldots, D_1(x_i), \ldots, x_n]=$$
$$D_{1}D_{2}([x_1, \ldots, x_n])-
\sum\limits_{i=2}^n \epsilon(D_1, D_2+X_i)\epsilon(D_2, X_i)[x_1,\ldots, D_2D_1(x_i), \ldots, x_n]$$

and
$$[D_{2}D_{1}(x_1), \ldots, x_n]=
D_{2}(D_{1}([x_1,\ldots, x_n])-
\sum\limits_{i=2}^n \epsilon(D_1,X_i)[x_1,\ldots, D_{1}(x_i),\ldots,x_n])=$$
$$D_{2}D_{1}([x_1, \ldots, x_n])-
\sum \epsilon(D_1,X_i)\epsilon(D_2,X_i) [x_1,\ldots, D_2D_{1}(x_i),\ldots,x_n].$$
Hence
$$[[D_{1},D_{2}](x_1), \ldots, x_n]=
D_{1}D_{2}([x_1, \ldots, x_n])- \epsilon(D_1,D_2) D_{2}D_{1}([x_1,\ldots, x_n])=[D_{1},D_{2}]([x_1, \ldots, x_n]).$$
 Similarly,
$$[[D_{1},D_{2}](x_1), \ldots, x_n]=\epsilon(D_1+D_2, X_i)[x_1, \ldots, [D_{1},D_{2}](x_i), \ldots, x_n].$$
Thus $[D_{1},D_{2}]\in{\rm C}(T)$ and $[{\rm Der}(T),{\rm
C}(T)]\subseteq {\rm C}(T).$

$(2)$ Similar to the proof of $(1)$.

$(3)$ Let $D_{1}\in{\rm C}(T),D_{2}\in{\rm Der}(T).$ For arbitrary
$x_1,\ldots, x_n \in T$ we have
$$D_{1}D_{2}[x_1, \ldots, x_n]=D_{1}\sum \epsilon (D_2, X_i) [x_1, \ldots, D_2(x_i), \ldots, x_n])=$$
$$\sum \epsilon (D_1+D_2, X_i) [x_1, \ldots, D_{1}D_2(x_i), \ldots, x_n]).$$
Therefore, $D_{1}D_{2}\in {\rm Der}(T).$

$(4)$ Let $D\in{\rm QC}(T).$ For arbitrary $x,y,z\in T,$ we have
$$[D(x_1), \ldots, x_n]= \epsilon(D,X_i)[x_1, \ldots, D(x_i),\ldots, x_n].$$
Hence
$$\sum  \epsilon(D,X_i)[x_1, \ldots, D(x_i), \ldots, x_n]=nD[x_1,\ldots, x_n].$$
Therefore $D\in {\rm QDer}(T)$ since $D^{'}=nD\in {\rm
C}(T)\subseteq {\rm End}(T)$.

$(5)$ Let $D_{1},D_{2} \in {\rm QC} (T).$ For arbitrary $x_1, \ldots, x_n \in T$
we have

$$[x_1, \ldots, [D_1,D_2](x_i),\ldots, x_n]=$$
$$\epsilon(X_i,D_1)\epsilon(X_i-x_1,D_2)[D_1(x_1),D_2(x_2), \ldots, x_n]-$$
$$\epsilon(D_1,D_2)\epsilon(X_i-x_1,D_2)\epsilon(X_i+D_2,D_1) [D_1(x_1),D_2(x_2), \ldots, x_n]=0.$$

Hence
$$\sum \epsilon(D_1+D_2, X_i) [x_1, \ldots, [D_1,D_2](x_i), \ldots, x_n]=0$$
and $[D_{1},D_{2}]\in {\rm QDer}(T)$.

$(6)$\quad Is obvious.\hfill$\Box$
 \vspace{0.3cm}


 \noindent{\bf Definition 2.3'.} {\it Let $T$ be a color $n$-ary algebra. Assume that $T$ satisfies the following property

$$ [x_1, \ldots, x_i, x_{i+1}, \ldots, x_n]= \gamma \epsilon(x_i,x_{i+1}) [x_1, \ldots, x_{i+1},x_i, \ldots, x_n], \forall i \in \{1, \ldots, n-1\}, \forall x_j \in T,$$
then\\
$(1)$ If $\gamma = 1$, we say that $T$ is a color $n$-ary commutative algebra,\\
$(2)$ If $\gamma=-1$, we say that $T$ is a color $n$-ary anticommutative algebra.}
\vspace{0.3cm}

\noindent{\bf Proposition 2.3.}
{\it Let $T$ be a color $n$-ary  (anti)commutative algebra. Then $GDer(T)=QDer(T)+QC(T).$}

\rm {\it Proof.}
It is easy to see that
$$[x_1, \ldots, x_i, \ldots, x_j, \ldots, x_n]=$$
$$ \gamma \epsilon(x_i+\ldots+x_{j-1}, x_j)\epsilon(x_i, x_{i+1}+ \ldots + x_{j-1})[x_1, \ldots, x_{i-1}, x_j, x_{i+1}, \ldots, x_{j-1}, x_{j}, x_{j+1}, \ldots, x_n].$$
Let $\epsilon_{i,j}= \gamma \epsilon(x_i +\ldots +x_{j-1},x_j)\epsilon(x_i, x_{i+1}+\ldots +x_{j-1}). $

Hence
$$D^{(n)}[x_1, \ldots, x_n]= \epsilon_{i,j} D^{(n)}[x_1, \ldots, x_j, \ldots, x_i,\ldots,  x_n]=$$
$$\epsilon_{i,j} \sum\limits_{t \leq i<j} \epsilon(D^{(t-1)}, X_t) [x_1, \ldots, D^{(t-1)}(x_t), \ldots, x_j, \ldots, x_i,\ldots, x_n]+$$
$$\epsilon_{i,j} \sum\limits_{i<t \leq j} \epsilon(D^{(t-1)}, X_t- x_i+x_j) [x_1, \ldots,x_i, \ldots, D^{(t-1)}(x_t), \ldots, x_i,\ldots, x_n]+$$
$$\epsilon_{i,j} \sum\limits_{i<j \leq t } \epsilon(D^{(t-1)}, X_t) [x_1, \ldots, D^{(t-1)}(x_t), \ldots, x_j, \ldots, x_i,\ldots, x_n]=$$

$$ \sum\limits_{t \neq i,j} \epsilon(D^{(t-1)}, X_t) [x_1, \ldots, D^{(t-1)}(x_t), \ldots, x_n]+$$
$$  \epsilon(D^{(j-1)}, X_i) [x_1, \ldots, D^{(j-1)}(x_i), \ldots,  x_j,\ldots, x_n]+$$
$$  \epsilon(D^{(i-1)}, X_j) [x_1, \ldots, x_i, \ldots, D^{(i-1)}(x_j), \ldots, x_n].$$

Hence, if $(D, D', \ldots, D^{(i-1)},  \ldots, D^{(n-1)},D^n)$ is a $(n+1)$-ary derivation,
then
$$(D, D', \ldots, \underbrace{D^{(i-1)}}\limits_{j}, , \ldots, \underbrace{D^{(j-1)}}\limits_{i}, \ldots, D^{(n-1)},D^{(n)})$$
and 
$$(D^{(i-1)}, D', \ldots, \underbrace{D^{(j-1)}}\limits_{i},  \ldots, \underbrace{D}\limits_{j} \ldots, D^{(n-1)},D^{(n)})$$
are  $(n+1)$-ary derivations.
It is easy to see that $D-D^{(i-1)} \in QC(T)$ and $nD - \sum\limits_{i=1}^{n} D^{(i-1)} \in QC(T).$

It is well known that every  permutation $\sigma \in S_{n-1}$ is a composition of transpositions $(i,j)$.
Hence we have
$$D^{(n)}[x_1, \ldots, x_n] = \sum\limits_{t=1}^{n} \epsilon(D, X_t) [x_1, \ldots, D^{(\sigma(t)-1)}(x_{t+1}), \ldots, x_n].$$

Obviously
$$\sum\limits_{i=1}^{n} \epsilon(D, X_i) [x_1, \ldots, (n-1)! (\sum\limits_{t=0}^{n-1} D^{(t)})(x_i),\ldots, x_n] = \sum\limits_{\sigma \in S_{n}} D^{(n)}[x_1, \ldots, x_n]= n! D^{(n)}[x_1, \ldots, x_n].$$

Therefore $\sum\limits_{i=1}^{n} D^{(i-1)} \in QDer(T).$
Now, every generalized derivation $D$ can be presented as a sum of a quasiderivation $( \sum\limits_{i=1}^{n} D^{(i-1)})/ n$
and an element of quasicentroid $D - (\sum\limits_{i=1}^{n} D^{(i-1)})/n.$
\hfill$\Box$\vspace{0.3cm}

\noindent{\bf Proposition 2.4.} {\it If $T$ is a color $n$-ary algebra,
then ${\rm QC}(T)+[{\rm QC}(T),{\rm QC}(T)]$ is a subalgebra of
${\rm GDer}(T)$.}

{\it Proof.}\quad By Lemma $2.2~(5)$ and $(6)$
we have
$${\rm QC}(T)+[{\rm QC}(T),{\rm QC}(T)]\subseteq {\rm GDer}(T)$$ and
\begin{eqnarray*}
&&\ \ [{\rm QC}(T)+[{\rm QC}(T),{\rm QC}(T)],{\rm QC}(T)+[{\rm QC}(T),{\rm QC}(T)]]\\
&&\subseteq[{\rm QC}(T)+{\rm QDer}(T),{\rm QC}(T)+[{\rm QC}(T),{\rm QC}(T)]]\\
&&\subseteq[{\rm QC}(T),{\rm QC}(T)]+[{\rm QC}(T),[{\rm QC}(T),{\rm
QC}(T)]]+[{\rm QDer}(T),{\rm QC}(T)]\\&&+[{\rm QDer}(T),[{\rm
QC}(T),{\rm QC}(T)]].\end{eqnarray*} Using the Jacobi identity, it is easy to verify that $[{\rm
QDer}(L),[{\rm QC}(L),{\rm QC}(L)]]\subseteq [{\rm QC}(L),{\rm
QC}(L)]$. Thus $$[{\rm
QC}(L)+[{\rm QC}(L),{\rm QC}(L)],{\rm QC}(L)+[{\rm QC}(L),{\rm
QC}(L)]]\subseteq {\rm QC}(L)+[{\rm QC}(L),{\rm QC}(L)].$$
\hfill$\Box$\vspace{0.3cm}
 \vspace{0.3cm}


\noindent{\bf Lemma 2.5.} {\it If $T$ is a color $n$-ary algebra,
then
$[{\rm C}(T),{\rm QC}(T)]\subseteq {\rm End}(T,{\rm Z}(T)).$
 Moreover, if ${\rm
Z}(T)=\{0\},$ then $[{\rm C}(T),{\rm QC}(T)]=\{0\}.$}

{\it Proof.}\quad Assume that
$D_{1}\in {\rm C}(T), D_{2}\in {\rm QC}(T)$ and for arbitrary $x_1, \ldots, x_n\in T$ we have
$$[[D_{1},D_{2}](x_1), \ldots, x_n]=[D_{1}D_{2}(x_1), \ldots, x_n ]- \epsilon(D_1,D_2) [D_{2}D_{1}(x_1), \ldots, x_n]=$$
$$D_{1}([D_{2}(x_1), \ldots, x_n])-\epsilon(D_2,x_1)[D_{1}(x_1),D_{2}(x_2), \ldots, x_n]$$
$$D_{1}([D_{2}(x_1), \ldots, x_n]-\epsilon(D_2,x_1)[x_1,D_{2}(x_2), \ldots, x_n])=0.$$
Hence
$[D_{1},D_{2}](x)\in {\rm Z}(T)$ and $[D_{1},D_{2}]\in {\rm
End}(T,{\rm Z}(T))$ as desired. Furthermore, if ${\rm Z}(T)=\{0\},$
it is clear that $[{\rm C}(T),{\rm QC}(T)]=\{0\}.$
\hfill$\Box$  \vspace{0.3cm}

\noindent{\bf  Definition 2.6.}
{\it Let $L$ be an $G$-graded algebra. If the relations
$$x\cdot y= \epsilon (x,y) y\cdot x, $$
$$\epsilon (z,x+w)  (((x\cdot y)\cdot w)\cdot z-(x\cdot y)\cdot(w\cdot z)))+$$
$$\epsilon (x,y+w)  (((y\cdot z)\cdot w)\cdot x-(y\cdot z)\cdot(w\cdot x)))+$$
$$\epsilon (y,z+w)  (((z\cdot x)\cdot w)\cdot y-(z\cdot x)\cdot(w\cdot y))=0,$$
hold in $L$ for all $x,y,z,w\in L,$ then we call $L$ a color Jordan algebra.}
\vspace{0.3cm}

\noindent{\bf  Proposition 2.7.} {\it Let $T$ be a color $n$-ary algebra,
then $End(T)$ with multiplication $$D_1 \bullet D_2 = D_1D_2 +\epsilon(D_1,D_2)D_2D_1$$
is a color Jordan algebra.} \vspace{0.3cm}

{\it Proof.}
Similar to Prop 2.9 in \cite{cml13}.\hfill$\Box$  \vspace{0.3cm}

\noindent{\bf  Corollary 2.8.} {\it Let $T$ be a color $n$-ary algebra. Then ${\rm QC}(T)$ endowed with the operation
$$D_{1}\bullet D_{2}=D_{1}D_{2}+\epsilon(D_1,D_2)D_{2}D_{1}~$$
is a color Jordan algebra.}

\rm {\it Proof.}\quad We need only to show that $D_{1}\bullet
D_{2}\in{\rm QC}(T)$. For arbitrary $D_{1},D_{2}\in{\rm QC}(T)$ we have

$$[D_{1}\bullet D_{2}(x_1), \ldots, x_n]=[D_{1}D_{2}(x_1), \ldots, x_n]+ \epsilon(D_1,D_2) [D_{2}D_{1}(x_1), \ldots, x_n]=$$
$$\epsilon(D_1,D_2+X_i)[D_{2}(x_1),\ldots, D_{1}(x_i), \ldots, x_n]+ \epsilon(D_2,X_i)[D_{1}(x_1), \ldots, D_{2}(x_i), \ldots, x_n]=$$
$$\epsilon(D_1,D_2+X_i)\epsilon(D_2,X_i)[x_1, \ldots, D_{2}D_{1}(x_i), \ldots, x_n]+ \epsilon(D_1+D_2,X_i)[x_1,\ldots, D_{1}D_{2}(x_i), \ldots, x_n]=$$
$$\epsilon(D_1+D_2,X_i)[x_1,\ldots, D_{1}\bullet D_{2}(x_i), \ldots,x_n].$$
Then $D_{1}\bullet D_{2}\in {\rm QC}(T)$ and ${\rm QC}(T)$ is a color Jordan algebra.
\hfill$\Box$  \vspace{0.3cm}

\noindent{\bf Theorem 2.9.}
{\it Let $T$ be a color $n$-ary algebra. Then the following statements hold:

$(1)$\quad ${\rm QC}(T)$ is a Lie algebra with $[D_{1},D_{2}]= D_{1}D_{2}-\epsilon(D_1,D_2)D_{2}D_{1}$ if and only if
${\rm QC}(T)$ is an associative algebra with respect to usual composition of operators;

$(2)$\quad If char $\mathbb {F}$ does not divide $n$ and ${\rm Z}(T)=\{0\}$, then
${\rm QC}(T)$ is a color Lie algebra if and only if
$~[{\rm QC}(T),{\rm QC}(T)]=0.$}

\rm {\it Proof.}\quad $(1)$ $(\Leftarrow)$ For arbitrary $D_{1},D_{2}\in
{\rm QC}(T)$ we have $D_{1}D_{2} \in {\rm QC}(T)$ and $D_{2}D_{1}\in
{\rm QC}(T)$, so $[D_{1},D_{2}]=D_{1}D_{2}-\epsilon(D_1,D_2)D_{2}D_{1}\in {\rm
QC}(T).$ Hence ${\rm QC}(T)$ is a Lie algebra.

$(\Rightarrow)$ Note that $D_{1}D_{2}=D_{1}\bullet D_{2}+\frac{[D_{1},D_{2}]}{2}$ and by Corollary $2.8$
$D_{1}\bullet D_{2}\in {\rm QC}(T), [D_{1},D_{2}]\in {\rm QC}(T)$.
It follows that $D_{1}D_{2}\in {\rm QC}(T)$ as desired.

$(2)$ $(\Rightarrow)$\quad  Let $D_{1},~D_{2}\in {\rm
QC}(T).$ Since ${\rm QC}(T)$ is a Lie algebra, for arbitrary $x_1, \ldots, x_n \in T$ we have
$$[[D_{1},D_{2}](x_1), \ldots, x_n]=\epsilon(D_1+D_2, X_i) [x_1,\ldots, [D_{1},D_{2}](x_i), \ldots, x_n].$$
From the proof of Lemma$~2.2~(5)$ it follows that
$$ \sum \epsilon(D_1+D_2, X_i)[x_1, \ldots, [D_{1},D_{2}](x_i), \ldots, x_n]=0.$$
Hence $n[[D_{1},D_{2}](x_1), \ldots, x_n]=0.$ Since char $\mathbb F$ does not divide $n$, we infer
$[[D_{1},D_{2}](x_1),\ldots, x_n]=0,$ i.e. $[D_{1},D_{2}]=0$.

$(\Leftarrow)$\quad Is clear.\hfill$\Box$


       \vspace{0.3cm}

\noindent {\bf \S 3\quad  Quasiderivations of color $n$-ary $\Omega$-algebras}

        \vspace{0.3cm}

 \noindent{\bf Definition 3.1.} {\it
For a (may be $n$-ary) multilinear polynomial $f(x_1, \ldots, x_n)$
we fix the order of indexes $\{i_1, \ldots, i_n\}$ of one non-associative word $[x_{i_1}\ldots x_{i_n}]_{\beta}$
from the polynomial $f$.
Here, $f = \sum\limits_{\beta, \sigma \in S_n} \alpha_{\sigma, \beta} [x_{\sigma(i_1)} \ldots x_{\sigma(i_n)}]_{\beta},$
where $\beta$ is an arrangement of brackets in the non-associative word.
For the  shift  $\mu_i: \{{j_1}, \ldots, {j_n} \} \mapsto \{j_1, \ldots j_{i+1}, j_i , \ldots,  j_n\}$
we define the element $\epsilon(x_{j_i},x_{j_{i+1}}).$
Now, for arbitrary non-associative word $ [x_{\sigma(i_1)} \ldots x_{\sigma(i_n)}]_{\beta}$ its order of indexes is a composition of suitable shifts $\mu_i$,
and for this word we set $\epsilon_{\sigma}$ defined as the product of corresponding $\epsilon(x_{j_i},x_{j_{i+1}}).$
Now, for the multilinear polynomial $f$, we define the color multilinear polynomial
$$f_{co}=\sum\limits_{\beta, \sigma \in S_n} \alpha_{\sigma, \beta} \epsilon_{\sigma} [x_{\sigma(i_1)} \ldots x_{\sigma(i_n)}]_{\beta}.$$}
\vspace{0.3cm}

 \noindent{\bf Definition 3.2.} {\it
Let $\Omega=\{ f_i\}$ be a family of $n$-ary multilinear polynimials.
A color  $n$-ary  $\Omega$-algebra $L$ is a color $n$-ary algebra satisfying the family of color  polinomials $\Omega_{co}=\left\{ (f_i)_{co} \right\}.$}
\vspace{0.3cm}

       In this section, we will prove that the quasiderivations of color  $n$-ary $\Omega$-algebra
       $T$ can be embedded as derivations in a larger  color $n$-ary $\Omega$-algebra and obtain
a direct sum decomposition of {\rm Der}($T$) when the center ${\rm Z}(T)$  of $T$ is equal to zero.
\vspace{0.3cm}

 \noindent{\bf Proposition 3.3.} {\it Let $T$ be a color $n$-ary $\Omega$-algebra over
       ${\mathbb F}$ and $t$ be an indeterminate. We define $\breve{T}:=
\{\Sigma(x\otimes t+y\otimes t^{n})| x,y\in T\}$.
Then $\breve{T}$
is a  color $n$-ary $\Omega$-algebra with the operation
$$[x_1\otimes t^{i_1},x_2 \otimes
t^{i_2}, \ldots ,x_n\otimes t^{i_n}]=[x_1,x_2, \ldots, x_n]\otimes t^{\sum i_j},$$
for $x_1,x_2, \ldots, x_n \in T, i_j\in\{1,n\}$}.

{\it Proof.}\quad Let the class of color $n$-ary $\Omega$-algebras be defined by family $\{ f_k \}$ of color multlilinear identities,
then for arbitrary $x_1, x_2, \ldots, x_m \in T$ and $i_j \in\{1,n\}$ we have

$$f_j(x_1\otimes t^{i_1},x_2\otimes t^{i_2}, \ldots, x_m\otimes t^{i_m})=
f_j(x_{1},x_{2}, \ldots, x_{m}) \otimes t^{\sum i_l}=0.$$

Therefore $\breve{T}$ is a color $n$-ary $\Omega$-algebra.

For the sake of convenience, we write $xt ~ (xt^{n})$ instead of $x\otimes t ~ (x\otimes t^{n}).$

If $U$ is a $G$-graded subspace of $T$ such that $T=U\oplus [T,\ldots,T],$ then
$$\breve{T}=Tt+Tt^{n}=Tt+Ut^{n}+[T,\ldots,T]t^{n}.$$
Now define a map $\varphi:{\rm QDer}(T)\rightarrow {\rm End}(\breve{T})$ by
$$\varphi(D)(at+ut^{n}+bt^{n})=D(a)t+D'(b)t^{n},$$
where $D\in {\rm QDer}($T$),$  $D'$ is a mapping related to $D$ by the definition of quasiderivation, $a\in
T,u\in U,b\in [T,\ldots,T]$. \vspace{0.3cm}

\noindent{\bf Proposition 3.4.} {\it Let $T,\breve{T},\varphi$ be as above. Then

$(1)$ $\varphi$ is injective and $\varphi(D)$ does not depend
on the choice of $D'$;

$(2)$ $\varphi({\rm QDer}(T))\subseteq {\rm Der}(\breve{T}).$}

{\it Proof.}\quad (1)\quad If $\varphi(D_{1})=\varphi(D_{2}),$ then
for all $a\in T,b\in [T,\ldots ,T]$ and $u\in U,$ we have
$$\varphi(D_{1})(at+ut^{n}+bt^{n})=\varphi(D_{2})(at+ut^{n}+bt^{n}),$$ or, in terms of $D_1, D_2,$
$$D_{1}(a)t+D'_{1}(b)t^{n}=D_{2}(a)t+D'_{2}(b)t^{n},$$ so $D_{1}(a)=D_{2}(a).$ Hence
$D_{1}=D_{2},$ and $\varphi$ is injective.

Suppose that there exists $D''$ such that
$$\varphi(D)(at+ut^{n}+bt^{n})=D(a)t+D''(b)t^{n},$$ and
$$\sum  \epsilon(D, X_i) [x_1, \ldots, D(x_i),\ldots, x_n]=D''([x_1, \ldots, x_n]),$$
then we have $$D'([x_1, \ldots, x_n])=D''([x_1, \ldots, x_n]),$$ thus $D'(b)=D''(b).$
Hence
$$\varphi(D)(at+ut^{n}+bt^{n})=D(a)t+D'(b)t^{n}=D(a)t+D''(b)t^{n},$$
which implies $\varphi(D)$ is determined only by $D$.

(2)\quad We have $[x_1t^{i_1}, \ldots, x_n t^{i_n}]=[x_1 , \ldots, x_n ]t^{\sum t_j}=0,$ for
all $\sum t_j \geq n+1$. Thus, to show $\varphi(D)\in {\rm
Der}(\breve{T}),$ we only need to check the validness of the
following equation:

$$\varphi(D)([x_1t, \ldots, x_nt])=\sum \epsilon(D, X_i) [x_1t,\ldots, \varphi(D)(x_it),\ldots, x_nt].$$

For arbitrary $x_1, \ldots, x_n \in T$ we have
$$\varphi(D)([x_1t,\ldots, x_it, \ldots, x_nt])=\varphi(D)([x_1,\ldots, x_n]t^n)=
D'([x_1, \ldots, x_n])t^n=$$
$$\sum \epsilon(D, X_i) [x_1, \ldots, D(x_i),\ldots, x_n]t^n=$$
$$\sum  \epsilon(D, X_i) [x_1t, \ldots, D(x_i)t,\ldots, x_nt]=$$
$$\sum \epsilon(D, X_i) [x_1t, \ldots, \varphi(D)(x_it),\ldots, x_nt].$$
Therefore, for all $D\in {\rm QDer}(T)$, we have $\varphi(D)\in {\rm
Der}(\breve{T})$.\hfill$\Box$ \vspace{0.3cm}

\

\noindent{\bf Proposition 3.5.} {\it Let $T$ be a $n$-ary $\Omega$-algebra such that
${\rm Z}(T)=\{0\}$ and $\breve{T},~\varphi$ be as defined above.
Then $${\rm Der}(\breve{T})=\varphi({\rm QDer}(T))\oplus {\rm
ZDer}(\breve{T}).$$}

{\it Proof.}\quad Since ${\rm Z}(T)=\{0\}$, we have ${\rm
Z}(\breve{T})=Tt^n.$
For all $g\in {\rm Der}(\breve{T}),$ we have
$g({\rm Z}(\breve{T}))\subseteq {\rm Z}(\breve{T}),$ hence
$g(Ut^n)\subseteq g({\rm Z}(\breve{T}))\subseteq {\rm
Z}(\breve{T})=Tt^n.$ Now define a map
$f:Tt+Ut^n+[T,\ldots,T]t^n\rightarrow Tt^n$ by
$$\ f(x)=\left\{\begin{array}{ll}g(x)\cap Tt^n,& x\in Tt ;\\
 g(x),& x\in Ut^n ;\\  0,& x\in [T,\ldots,T]t^n.\end{array}\right.$$

It is clear that $f$ is linear. Note that
$$f([\breve{T}, \ldots, \breve{T}])=f([T,\ldots ,T]t^n)=0,$$
$$~[\breve{T}, \ldots, f(\breve{T}), \ldots, ,\breve{T}]
\subseteq [Tt+Tt^n,\ldots, Tt^n,\ldots, Tt+Tt^n]=0,$$
hence $f\in {\rm ZDer}(\breve{T}).$
Since $$(g-f)(Tt)=g(Tt)-g(Tt)\cap Tt^{n}=g(Tt)-Tt^{n}\subseteq Tt,~
(g-f)(Ut^n)=0,$$ and
$$(g-f)([T, \ldots, T]t^n)=g([\breve{T},\ldots ,\breve{T}])\subseteq
[\breve{T},\ldots, \breve{T}]=[T, \ldots, T]t^n,$$ there exist $D,~D'\in
{\rm End}(T)$ such that for all $a\in T,~b\in [T, \ldots, T]$,
$$(g-f)(at)=D(a)t,~ (g-f)(bt^n)=D'(b)t^n.$$ Since $(g-f)\in {\rm
Der}(\breve{T})$ and by the definition of ${\rm Der}(\breve{T})$, we
have
$$\sum \epsilon(g-f, A_i) [a_1, \ldots, (g-f)(a_it), \ldots, a_nt]=(g-f)([a_1t, \ldots, a_nt]),$$
for all $a_1, \ldots, a_n\in T.$
Hence
$$\sum \epsilon(D, A_i) [a_1, \ldots, D(a_i), \ldots,a_n]=D'([a_1,  \ldots, a_n]).$$ Thus
$D\in {\rm QDer}(T).$ Therefore, $g-f=\varphi(D)\in \varphi({\rm
QDer}(T))$, so ${\rm Der}(\breve{T})\subseteq \varphi({\rm
QDer}(T))+{\rm ZDer}(\breve{T}).$ By Proposition $3.4~(2)$ we have
${\rm Der}(\breve{T})=\varphi({\rm QDer}(T))+{\rm ZDer}(\breve{T}).$

For any $f\in \varphi({\rm QDer}(T))\cap{\rm ZDer}(\breve{T})$
there exists an element $D\in {\rm QDer}(T)$ such that
$f=\varphi(D).$ Then
$$f(at+ut^n+bt^n)=\varphi(D)(at+ut^n+bt^n)=D(a)t+D'(b)t^n,$$ where $a\in T,b\in [T, \ldots, T].$

On the other hand, since $f\in {\rm ZDer}(\breve{T}),$ we have
$$f(at+bt^n+ut^n)\in {\rm Z}(\breve{T})=Tt^n.$$ That is to say,
$D(a)=0,$ for all $a\in T$ and so $D=0.$ Hence $f=0.$

Therefore ${\rm Der}(\breve{T})=\varphi({\rm QDer}(T))\oplus {\rm
ZDer}(\breve{T})$ as desired. \hfill$\Box$\vspace{0.3cm}



\medskip
\noindent {\bf \S 4\quad  Algebras with the condition $End=QDer.$}
        \vspace{0.3cm}

\medskip

Throughout the section all spaces of algebras are assumed finite-dimensional \textbf{over a field of characteristic 0}.
An $n$-ary algebra $L$ with multiplication $[\cdot,\ldots,\cdot]$ will be called {\it commutative} if it satisfies the identity
$$[x_{\sigma(1)},x_{\sigma(2)},\ldots,x_{\sigma(n)}] = [x_1,x_2,\ldots,x_n]$$
for any $\sigma \in S_n$, and {\it anticommutative}, if it satisfies the identity
$$[x_{\sigma(1)},x_{\sigma(2)},\ldots,x_{\sigma(n)}] = (-1)^{sgn(\sigma)}[x_1,x_2,\ldots,x_n]$$
for any $\sigma \in S_n.$

For any commutative $n$-ary algebra $L$ we may consider its multiplication as a linear mapping $\mu:S^n(L) \rightarrow L$, where $S^n(L)$ is the {\it $n$th symmetric power} of $L$, which is the space of elements of the symmetric algebra $Sym(L)$, homogeneous of degree $n$. Analogously, for any anticommutative $n$-ary algebra $L$, we may consider its multiplication as a linear mapping $\mu: \bigwedge^n(L) \rightarrow L,$ where $\bigwedge^n(L)$ is the {\it $n$th exterior power} of $L$, which is the space of elements of Grassmann algebra $\bigwedge(L)$, homogeneous of degree $n$.

An $n$-ary algebra $L$ with anticommutative multiplication $[\cdot,\ldots,\cdot]$ is called a {\it Filippov (or $n$-Lie) algebra} if it satisfies the following identity:
$$[[x_1,\dots,x_n],y_2,\dots,y_n] = \sum_{i=1}^n[x_1, \dots, [x_i,y_2,\dots,y_n],\dots,x_n].$$

Filippov algebras were introduced in \cite{Fil}.
In his article, Filippov introduced the notion of a $n$-Lie algebra and proved some properties of them. He also obtained the classification of anticommutative $n$-ary algebras of dimension $n$ and $n+1$.
Here we list some important $n$-Lie and anticommutative algebras of dimension $\leq n+1$:

$1)$ Up to isomorphism, there is only one $n$-ary anticommutative algebra $A_n$ in dimension $n.$ On the basis elements $e_1, \dots, e_n$ of $A_n$ we define the product in the following way:
$$[e_1,\dots,e_n] = e_1.$$
Then $A_n$ is an anticommutative (even Filippov) algebra and any $n$-dimensional $n$-ary anticommutative algebra is isomorphic to $A_n.$

$2)$ Up to isomorphism, there is only one perfect $n$-Lie algebra $D_{n+1}$ of dimension $n+1$.
On the basis elements $e_1, \dots, e_{n+1}$ of $D_{n+1}$ we define the product in the following way:
$$ [e_1, \dots, \hat{e_i}, \dots, e_{n+1}] = (-1)^{n+i+1}e_i.$$

$3)$ Let $L$ be an $n$-ary $n+1$-dimensional algebra with basis $e_1, \dots, e_n.$ Let
\begin{equation}
\label{ei}
e^i = (-1)^{n+i+1}[e_1,\dots,\hat{e_i},\dots,e_n], i = 1, \dots, n+1.
\end{equation}
Then, multiplication in $L$ is defined by the matrix $B = (\beta_{ij})$ which is given by $e^i = \beta_{1i}e_1 + \dots + \beta_{n+1 i}e_{n+1}$, or in terms of matrices as
\begin{eqnarray}
\label{L_B}
(e^1,\ldots,e^{n+1}) = (e_1,\ldots,e_{n+1})B.
\end{eqnarray}

It is obvious that the rank of $B$ is equal to the dimension of $L^2 = [L,\ldots,L].$ 
An $(n+1)$-dimensional anticommutative algebra with multiplication defined by (\ref{ei}) and (\ref{L_B}) will be denoted by $L_B.$ It is easy to see that $D_{n+1}$ is in fact an algebra $L_I$, where $I$ is an identity matrix of order $n+1.$

Also later in the discussion we will need the description of 1-dimensional $n$-ary algebras. One can easily see that the multiplication in such algebra is completely defined by an element $\alpha \in F$, where $F$ is the base field, for it is enough to determine
$$[v,v,\ldots,v] = \alpha v$$
for any nonzero $v \in L$ and extend the multiplication linearly. We denote such algebras by $L_{\alpha}$. It is also easy to see that $L_{\alpha} \cong L_{\beta}$ for $\alpha, \beta \neq 0$ if and only if the polynomial $x^{n-1} - \frac{\alpha}{\beta}$ has a root in $F.$

\

{\bf Definition 4.1.} Let $A$ be an $n$-ary algebra.
A pair of linear mappings $(d, f)$ is called a {\it quasiderivation} of $A$ if for arbitrary $a_1, \dots, a_n \in A$ we have
\begin{equation}
\label{QDer_def}
d([a_1,\dots,a_n]) = \sum_{i=1}^n[a_1, \dots, a_{i-1}, f(a_i), a_{i+1}, \dots, a_n].
\end{equation}
The image of projection of $QDer(A)$ onto the first coordinate will be denoted as $QDer_{KS}(A)$ and will be called a \emph{space of quasiderivations in a sense of Kaygorodov and Shestakov} \cite{kay11aa,lesha12},
and the image of projection of $QDer(A)$ onto the second coordinate will be denoted as $QDer_{LL}(A)$ and will be called a \emph{space of quasiderivations in a sense of Leger and Luks} \cite{LL00}.

In \cite{LL00} Leger and Luks completely described all binary algebras with the property $QDer_{LL}(A) = End(A).$ In this section we  describe all (anti)commutative $n$-ary algebras such that any its endomorphism is a quasiderivation in a sense of LL or in a sence of KS.

Let $(L, \mu)$ be a commutative algebra. For any $\phi \in End(L)$ we define $\phi^{*} \in End(Sym^n(L))$ by
$$\phi^{*}(x_1 \cdot \ldots \cdot x_n) = \sum_{i=1}^n x_1 \cdot \ldots \cdot x_{i-1} \cdot \phi(x_i) \cdot x_{i+1} \cdot \ldots \cdot x_n.$$
Analogously, let $(L, \mu)$ be an anticommutative algebra. For any $\phi \in End(L)$ we define $\phi^{*} \in End(\bigwedge^n(L))$ by
$$\phi^{*}(x_1\wedge \ldots \wedge x_n) = \sum_{i=1}^n x_1 \wedge \ldots \wedge x_{i-1} \wedge \phi(x_i) \wedge x_{i+1} \wedge \ldots \wedge x_n.$$

\

{\bf Lemma 4.2.} 
{\it Let $L$ be a commutative or anticommutative algebra. A mapping $f \in End(L)$ is a quasiderivation (in a sence of Leger and Luks) of $L$ if and only if $$f^{*}(\ker(\mu)) \subseteq \ker(\mu).$$}

{\it Proof.} It is evident that $f$ is a quasiderivation (in a sence of Leger-Luks) of $L$ if and only if there is a $d \in End(L)$ such that $f^{*}\mu = \mu d.$ By the Homomorphism theorem, this is equivalent to $\ker(\mu) \subseteq \ker(f^{*}\mu),$ and this condition is easily seen to be equivalent to $f^{*}(ker(\mu)) \subseteq ker(\mu).$ \hfill$\Box$ \vspace{0.3cm}

Now, let $L$ be an anticommutative $n$-ary algebra such that $QDer_{LL}(L) = End(L).$ Then, by previous lemma, $\ker(\mu)$ is a submodule of $\bigwedge^n(L)$ with respect to the action
\begin{equation}
\label{action}
\phi\cdot(x_1\wedge \ldots \wedge x_n) = \phi^{*}(x_1 \wedge \ldots \wedge x_n).
\end{equation}
The next lemma characterizes all submodules of $\bigwedge^n(L):$ with respect to this action:

\

{\bf Lemma 4.3.} 
{\it $\bigwedge^n(L)$ is an irreducible $End(L)$-module with respect to action (\ref{action}).}

{\it Proof.} Let $e_1, \dots, e_m, m \geq n$ be a basis for $L$ (if $m < n$ we have $\bigwedge^n(L) = 0$, so in this case the statement is trivial). We describe the basis for $\bigwedge^n(L).$ Let $I$ be any ordered $n$-tuple. that is, $I = \{i_1, \dots, i_n\}, i_1 < \ldots < i_n.$ Then it is clear that the elements $e_I = e_{i_1}\wedge \ldots \wedge e_{i_n}$ form the basis of $\bigwedge^n(L).$ Let $e_{ij}$ be standard matrix units in the basis $\{e_i\}$, that is, $e_{ij}e_k = \delta_{ik}e_j.$ Next we ask ourself the following question: how do $e_{ij}^{*}$ act on the basis elements $e_{I}?$\\
Let at first $i \neq j$. Then it is easy to see that
\begin{equation}
\label{act_ij}
e_{ij}^{*}(e_I) =
\begin{cases}
0, \text{if } i \notin I,\\
0, \text{if } i, j \in I,\\
\pm e_{I\setminus \{i\} \cup \{j\}}, \text{if } i \in I, j \notin I.
\end{cases}
\end{equation}

In the case when $i = j$ we have
\begin{equation}
\label{act_ii}
e_{ii}^{*}(e_I)=
\begin{cases}
0, \text{if } i \notin I,\\
e_I, \text{if } i \in I.
\end{cases}
\end{equation}

Now, let $0 \neq M \subseteq \bigwedge^n(A)$ be a submodule with respect to action (\ref{action}), and let $0 \neq x = \sum \alpha_I e_I \in M.$
Take an arbitrary $I_0= ( i_1, \ldots, i_n )$
such that $\alpha_{I_0} \neq 0.$
Then, by (\ref{act_ii}) we have
$$\frac{1}{\alpha_I}e_{i_1i_1}^{*}(  \ldots e_{i_ni_n}^{*} (x) \ldots ) = e_{I_0} \in M.$$

Now, using (\ref{act_ij}) it easy to see that by iteratively applying $\pm e_{ij}^{*}$ to $e_{I_0}$, we can show that $e_{J} \in M$ for any $n$-element ordered subset $J \subseteq \{e_1, \dots e_m\}.$ Therefore, $M = \bigwedge^n(L).$ \hfill$\Box$ \vspace{0.3cm}

Now let $L$ be a commutative $n$-ary algebra such that $QDer_{LL}(L) = End(L).$ Then again $\ker(\mu)$ is a submodule of $Sym^n(L)$ with respect to the action
\begin{equation}
\label{action_sym}
\phi\cdot(x_1 \cdot \ldots \cdot x_n) = \phi^{*}(x_1 \cdot \ldots \cdot x_n).
\end{equation}

The next lemma characterizes all submodules of $Sym^n(L)$ with respect to this action:

\

{\bf Lemma 4.4.}
{\it $Sym^n(L)$ is an irreducible $End(L)$-module with respect to action (\ref{action_sym}).}

{\it Proof.} Let $e_1, \dots, e_m$ be a basis for $L$. We describe the basis for $Sym^n(L).$ Let $I =(i_1, \ldots, i_m)$ be an $m$-tuple of natural numbers such that $i_1 + \ldots + i_m = n.$ Since the algebra $Sym(L)$ is just the algebra of polynomials in commuting variables $e_1, \ldots, e_m$, it is clear that the elements $e_I = e_1^{i_1} \cdot \ldots \cdot e_m^{i_m}$ form the basis of $Sym^n(L).$ Let $e_{ij}$ be standard matrix units in the basis $\{e_i\}$, that is, $e_{ij}e_k = \delta_{ik}e_j.$
Let us fix a basis element $e_I = e_1^{k_1} \cdot \ldots \cdot e_m^{k_m}.$
Then it is easy to see that
\begin{equation}
\label{act_ij_sym}
e_{ij}^{*}(e_I) = k_i e_{I'}, \text{ where } I' = (k_1, \dots, k_i -1, \dots, k_j +1, \dots, k_m), i, j = 1, \dots, m.
\end{equation}
Now, let $0 \neq M $ be a submodule of $Sym^n(A)$ with respect to the action (\ref{action_sym}).

Let $0 \neq x = \sum \alpha_I e_I \in M$.
Let $I_0=(i_1, \ldots, i_m)$ be such that $\alpha_{I_0} \neq 0$ and that the $m$th component $i_m$ is the least of all $m$ths components of all $I$ such that $\alpha_I \neq 0.$ 
Now,
$$\frac{1}{\alpha_{I_0} i_1! \ldots i_{m-1}!}(e_{1m}^*)^{i_1}( \ldots (e_{m-1,m}^*)^{i_{m-1}} (x) \ldots) = e_{(0,\ldots,0,n)} \in M.$$
Finally, using (\ref{action_sym}) and the fact that $e_{(0,\ldots,0,n)} \in M$ we can easily see that $e_J \in M$ for any $J,$ and $M=Sym^n(L).$

\hfill$\Box$

\

Now we can bound the dimension of $L:$

\

{\bf Lemma 4.5.} {\it Let $(L, \mu)$ be an $n$-ary algebra with $QDer_{LL}(L) = End(L).$

$1)$ If $L$ is anticommutative, then either $\mu = 0$ or $\dim(L) \leq n+1;$

$2)$ If $L$ is commutative, then either $\mu = 0$ or $\dim(L) = 1.$}

{\it Proof.}
$1)$ By Rank-nullity theorem, we have $\dim(L) \geq \dim(\bigwedge^n(L)) - \dim(\ker(\mu)).$ By lemmas 4.2 and 4.3  the only options we have are $\ker(\mu) = \bigwedge^n(L)$ (that is, $\mu = 0$) or $\ker(\mu) = 0.$ Suppose that the second option holds. Then we must have
\begin{equation}
\label{dim_ineq}
\dim(L) \geq \dim(\wedge^n(L)) = {\dim (L) \choose n}.
\end{equation}
This inequality can only hold when $\dim(L) \leq n+1.$

$2)$ Again, by Rank-nullity theorem, if $\mu \neq 0$, we infer
\begin{equation}
\label{dim_ineq_comm}
\dim(L) \geq \dim(Sym^n(L)) = {\dim(L) + n - 1 \choose n}.
\end{equation}
Let $\dim(L) = m > 1.$ Then it is easy to see that
$${m + n - 1 \choose n} = \frac{m\cdot(m+1)\cdot(m+2)\cdot\ldots\cdot(m+n-1)}{n!} > $$
$$> \frac{m(1+1)\cdot(1+2)\cdot\ldots\cdot(1+n-1)}{n!} = m.$$
Therefore, the inequality (\ref{dim_ineq_comm}) can only hold when $m = \dim(L) =  1.$
\hfill$\Box$

\

Now we can prove our main theorem.

\


{\bf Theorem 4.6.} {\it Let $L$ be an $n$-ary algebra such that $QDer_{LL}(L) = End(L)$. Then

 $1)$ If $L$ is commutative, then either $L$ has zero product, or $L \cong L_{\alpha}$ for some $\alpha$ in base field $F.$ Moreover, if the base field $F$ of $L$ is algebraically closed, then one can define a binary multiplication $\cdot$ on $L$ such that $(L,\cdot) \cong F$ and $[x_1,\ldots,x_n] = x_1\cdot\ldots\cdot x_n$, where $[\cdot,\ldots,\cdot]$ is a multiplication in $L$ and  $x_1, \dots, x_n \in L.$

 $2)$ If $L$ is anticommutative, then either $L$ has zero product, or $L$ is isomorphic to either $A_n$ or $L_B$ for nondegenerate $(n+1) \times (n+1)$-matrix $B$. Moreover, in the last case if $L$ is a Filippov algebra and the base field is algebraically closed, then} $L \cong D_{n+1}.$}

{\it Proof.} $1)$ Let $\mu \neq 0.$ By the previous lemma and the classification of 1-dimensional $n$-ary algebras, we infer $L \cong L_{\alpha}$ for some $ 0 \neq \alpha \in F.$ Since $L$ has only scalar linear mappings, it is obvious that $QDer_{LL}(L) = End(L).$ Now, if $F$ is algebraically closed, by remark above $L \cong L_1$ and there exists $e \in L$ such that $[e,\ldots,e] = e.$ Now we can define binary multiplication $\cdot$ in $F$ by setting $e\cdot e = e$ and extending linearly. It is obvious that $[x_1,\ldots,x_n] = x_1 \cdot \ldots \cdot x_n$.

$2)$ Suppose that $L$ has nonzero product. Using (\ref{dim_ineq}), we see that we have only two opportunities for $L:$

$2.1.$ $\dim(L) = n, \dim L^2 = 1.$ By Filippov's result, $L \cong A_n.$ Since $\ker(\mu)$ is zero, by lemma 4.2 we have $QDer_{LL}(A_n) = End(A_n).$

$2.2.$ $\dim(L) = n+1 = \dim(Im(\mu)) = dim(L^2).$ Therefore, $L$ is isomorphic to an algebra $L_B$ for some $(n+1) \times (n+1)$-matrix $B$ such that $rank(B) = dim(L^2) = n+1$. The Rank-nullity theorem implies that the kernel of $\mu$ is zero, therefore by lemma 4.2 $QDer_{LL}(L) = End(L).$ Now suppose that $\dim L = n+1,$ $L$ is a Filippov algebra and base field is algebraically closed. Checking with the list of all $(n+1)$-dimensional Filippov algebras over an algebraically closed field, we only have one possibility: $L \cong D_{n+1}$. 
\hfill$\Box$ \vspace{0.3cm}

Now we characterize (anti)commutative $n$-ary algebras with $QDer_{KS}(L) = End(L).$
First, we prove the obvious lemma:

\

{\bf Lemma 4.7.} {\it Let $A$ be an $n$-ary algebra such that $QDer_{KS}(A) = End(A).$ Then either $A^2 = 0$ or $A^2 = A$.}

{\it Proof.} From (\ref{QDer_def}) it follows that $A^2$ is invariant under all $QDer_{KS}(A) = End(A).$ But the only two $End(A)$-invariant spaces for any $A$ is $A$ itself and zero subspace. \hfill$\Box$ \vspace{0.3cm}

Now we reduce the problem to the previous one:

\

{\bf Lemma 4.8.} {\it Let $A$ be a  $n$-ary algebra, then

$1)$ Let $(d_1, f_1) \dots, (d_k, f_k) \in QDer(A).$ Suppose that $f_1, \dots, f_k$ are linearly dependent. Then the restrictions of $d_1, \ldots, d_k$ on $A^2$ are linearly dependent.

$2)$ Let $A$ be an $n$-ary algebra such that $A = A^2.$ Then $\dim QDer_{KS}(A) \leq \dim QDer_{LL}(A).$

$3)$ Let $A$ be an $n$-ary algebra with nonzero multiplication such that $QDer_{KS}(A) = End(A).$ Then $QDer_{LL}(A) = End(A).$}

{\it Proof.} $1)$ Suppose that $\sum_{i=1}^k \alpha_i f_i = 0,$ where $\alpha_i \in F.$ Then
$$(\sum_{i=1}^k \alpha_i d_i)[x_1, \dots, x_n] = \sum_{i=1}^n[x_1, \dots, x_{i-1}, (\sum_{i=1}^k \alpha_i f_i)(x_i), x_{i+1}, \dots, x_n] = 0,$$
which means that the restrictions of $d_1, \dots d_n$ on $A^2$ are linearly dependent.

$2)$ Let $d_1, \dots,d_k$ be a basis for $QDer_{KS}(A)$ and let $f_1, \dots, f_k$ be the corresponding Leger-Luks quasiderivations. Then from the previous point it follows that $f_1, \dots, f_k$ are linearly independent.

$3)$ The statement easily follows from the lemma 4.7 and the point $2).$ \hfill$\Box$ \vspace{0.3cm}

Now we can describe $n$-ary (anti)commutative algebras $L$ such that $QDer_{KS}(L) = End(L).$

\

{\bf Theorem 4.9.}
{\it Let $L$ be an $n$-ary algebra such that $QDer_{KS} (L) = End(L).$ Then

$1)$ If $L$ is commutative, then $L \cong L_{\alpha}$ for some $\alpha \in F.$ Moreover, if the base field $F$ of $L$ is algebraically closed, then one can define a binary multiplication $\cdot$ in $L$ such that $(L,\cdot) \cong F$ and $[x_1,\ldots,x_n] = x_1\cdot\ldots\cdot x_n$, where $[\cdot,\ldots,\cdot]$ is a multiplication in $L$ and  $x_1, \dots, x_n \in L.$

$2)$ If $L$ is anticommutative, then $L \cong L_B$ for nondegenerate $(n+1) \times (n+1)$-matrix $B$  or $L$ has zero multiplication. Moreover, if $L$ is a Filippov algebra and the base field is algebraically closed, then $L \cong D_{n+1}.$}

{\it Proof.} $1)$ Obviously follows from theorem 4.6.\\
 $2)$ Suppose that $L$ has nonzero multiplication. Lemma 4.8 and theorem 4.6  imply that $L$ is isomorphic either to $A_n$, or to $L_B$ for nondegenerate $B.$ Since $\dim A_{n}^2 = 1,$ lemma 4.7 implies that $QDer_{KS}(A_{n}) \neq End(A_n).$
Consider the case $L \cong L_B,$ where $B$ is a nondegenerate matrix of order $n+1.$ Let $e_{ij}$ and $e^{ij}, i, j = 1, \dots, n$ denote the matrix units in bases $e_1, \dots, e_n$ and $e^1, \dots, e^n$ correspondingly, that is, $e_{ij}(e_k) = \delta_{ik}e_j, e^{ij}(e^k) = \delta_{ik}e^j.$ To prove that $QDer_{KS}(L_B) = End(L_B)$ it suffices to check that $e^{ij} \in QDer_{KS}(L_B), i, j = 1, \dots, n.$
One can check directly that the mapping $e^{ij}$ is a quasiderivation in a sense of Kaygorodov-Shestakov with the corresponding Leger-Luks quasiderivation $-e_{ji}$. Indeed, for $i \neq j$ we have

$$e^{ij}  [e_1, \dots, \hat{e_i}, \dots, e_n] = (-1)^{n+i+1} e^{ij}(e^i) = (-1)^{n+i+1} e^j =$$
$$(-1)^{i - j}[e_1, \dots, \hat{e_j}, \dots, e_{n+1}]=$$
$$[e_1, \dots, e_{j-1}, e_i, e_{j+1}, \dots, \hat{e_i}, \dots, e_n] = $$
$$\sum_{k = 1}^n [e_1, \dots, -e_{ji}(e_k), \dots, \hat{e_i}, \dots, e_n]$$
and
$$e^{ij}[e_1, \dots, \hat{e_l}, \dots, e_n]=(-1)^{n-l+1}e^{ij}(e^l)= 0 =$$
$$\sum_{k = 1}^n [e_1, \dots, -e_{ji}(e_k), \dots, \hat{e_l}, \dots, e_n].$$

It is easy to see that these equalities remain correct in the case $i = j.$
Therefore, $QDer_{KS}(L_B) = End(L_B).$
\hfill$\Box$ \vspace{0.3cm}

\

\noindent {\bf \S 5\quad Appendix A.} 
        \vspace{0.3cm}

\medskip

In the end of the paper, we remark that it would be interesting to describe all finite-dimensional $n$-ary algebras with the property $QDer = End.$ The attempt of doing so in the case of Lie triple systems ($n$ = 3) was made in the paper \cite{cml15}. Unfortunately, authors made a mistake in the classification of the irreducible submodules of $T \otimes T \otimes T$ with respect to the action $f \cdot x \otimes y \otimes z = f^{*}(x \otimes y \otimes z)$, where $T$ is a Lie triple system, and the action of $End(T)$ on $T \otimes T \otimes T$ is defined similar to the actions (\ref{action}) and (\ref{action_sym}). Particularly, it is proposed that the sets 
\begin{eqnarray}\label{module}
(T \otimes T \otimes T)^{+} = span(x \otimes y \otimes z + y \otimes x \otimes z : x, y, z \in T)
\end{eqnarray}
and 
\begin{eqnarray}\label{module2}(T \otimes T \otimes T)^{-} = span(x \otimes y \otimes z - y \otimes x \otimes z : x, y, z \in T)
\end{eqnarray}
are the only two $End(T)$-submodules with respect to the above action. However, it is easy to see that the space 
$$span(x \otimes y \otimes z + x \otimes z \otimes y : x, y, z \in T )$$ 
is an $End(T)$-submodule with respect to the action above and that it does not coincide with any of the two spaces above.
In addition, 
the submodule (\ref{module}) has a non-trivial $End(T)$-submodule:
\begin{eqnarray*}(T \otimes T \otimes T)^{*} = span( \sum_{\sigma \in S_3} x_{\sigma(1)} \otimes x_{\sigma(2)} \otimes x_{\sigma(3)} : x_1, x_2, x_3 \in T);
\end{eqnarray*}
the submodule (\ref{module2}) has a non-trivial $End(T)$-submodule:
\begin{eqnarray*}(T \otimes T \otimes T)^{**} = span( \sum_{\sigma \in S_3} (-1)^{\sigma}x_{\sigma(1)} \otimes x_{\sigma(2)} \otimes x_{\sigma(3)} : x_1, x_2, x_3 \in T).\end{eqnarray*}

We also remark that that the proof of irreducibility of $Sym^n(L)$ and $\bigwedge^n(L)$ as $End(L)$-modules carries over verbatim to the case of $n$-ary color algebras.

       \vspace{0.3cm}



\begin{thebibliography}{99}


\bibitem{poj08patricia}  {\it Beites P., Pozhidaev A.,} On simple Filippov superalgebras of type A(n,n),
Asian-Eur. J. Math., 1 (2008), 4, 469--487.



\bibitem{ck10}  {\it Cantarini N., Kac V.,}
Classification of simple linearly compact $n$-Lie superalgebras, Comm. Math. Phys., 298 (2010), 3, 833--853.


\bibitem{cml13} {\it Chen L.,  Ma Y., Ni L.,}
Generalized Derivations of Lie color algebras, Results Math., 63 (2013),  3-4, 923--936.

\bibitem{cml14} {\it Chen L.,  Ma Y., Zhou J.,}
Generalized Derivations of Hom-Lie Superalgebras,
arXiv:1406.1578

\bibitem{cml15} {\it Chen L.,  Ma Y., Zhou J.,}
Generalized Derivations of Lie triple systems,
arXiv:1412.7804


\bibitem{dzhuma}
{\it Dzhumadilʹdaev A.,}
Cohomologies of colour Leibniz algebras: pre-simplicial approach,
Lie theory and its applications in physics, III (Clausthal, 1999), 124–136, World Sci. Publ., River Edge, NJ, 2000.



\bibitem{Fil}
{\it Filippov V.,}
$n$-Lie algebras,
Siberian Mathematical Journal, 26 (1985), 6, 126--140.


\bibitem{fil1} {\it Filippov V.,}
 $\delta$-Derivations of  Lie algebras,
Siberian Mathematical Journal, 39 (1998), 6, 1218--1230.

\bibitem{fil2} {\it Filippov V.,}
 $\delta$-Derivations of prime Lie algebras,
Siberian Mathematical Journal, 40 (1999), 1, 174--184.

\bibitem{fil3} {\it Filippov V.,}
$\delta$-Derivations of prime alternative and Malcev algebras,
Algebra and Logic, 39 (2000), 5, 354--358.

\bibitem{GP03} {\it Jimenez-Gestal C.,  Perez-Izquierdo J.M.,}
Ternary derivations of generalized Cayley–Dickson algebras, Comm. Algebra, 31 (2003), 10, 5071--5094.


\bibitem{GP08} {\it Jimenez-Gestal C.,  Perez-Izquierdo J.M.,}
Ternary derivations of finite-dimensional real division algebras, Lin. Alg. Appl., 428 (2008), 8-9, 2192--2219.

\bibitem{her}
{\it Herstein I.}, Jordan derivations of prime rings, Proc. Amer. Math. Soc., 8 (1957), 1104--1110.

\bibitem{kay1}  {\it Kaygorodov I.,} $\delta$-Derivations of simple finite-dimensional Jordan superalgebras,
Algebra and Logic, 46 (2007), 5, 318--329.


\bibitem{kay_lie} {\it Kaygorodov I.,}
$\delta$-Derivations of classical Lie superalgebras,
Sib. Math. J., 50 (2009), 3, 434-449.


\bibitem{kay_lie2}  {\it Kaygorodov I.,}
$\delta$-Superderivations of simple finite-dimensional Jordan and Lie superalgebras,
Algebra and Logic,  49 (2010), 2, 130--144.



\bibitem{kay12mz}   {\it Kaygorodov I.,}
$\delta$-Superderivations of semisimple finite-dimensional Jordan  superalgebras,
Mathematical Notes, 91 (2012), 2, 187--197.

\bibitem{kay12izv}   {\it Kaygorodov I.,}
$\delta$-Derivations of $n$-ary algebras,
Izvestiya: Mathematics, 76 (2012), 5, 1150--1162.




\bibitem{kay11aa} {\it Kaygorodov I.,}
$(n+1)$-Ary derivations of simple $n$-ary algebras, Algebra and Logic, 50 (2011), 5, 470--471.

\bibitem{kay14sp} {\it Kaygorodov I.,}
$(n+1)$-Ary derivations of simple Malcev algebras,
St. Peterburg Math. Journal, 23 (2014), 4, 575--585.


\bibitem{kay14mz}
{\it Kaygorodov I.,}
$(n+1)$-Ary derivations of semisimple Filippov algebras,
Math. Notes, 96 (2014), 2, 208--216.


\bibitem{KP}
{\it Kaygorodov I., Popov Yu.,}
Alternative algebras admitting derivations with invertible values and invertible derivations,
Izv. Math., 78 (2014), 5, 922--935.



\bibitem{KN03} {\it Komatsu H., Nakajima A.,}
 Generalized derivations of associative algebras, Quaest. Math., 26 (2003), 2, 213--235.






\bibitem{LL00}
{\it Leger G., Luks E.,}
Generalized derivations of Lie algebras,
J. Algebra, 228 (2000), 1, 165--203.




\bibitem{mikh85} {\it Mikhalev A.,} Subalgebras of free colored Lie superalgebras, Mathematical Notes, 37 (1985), 5, 356--360.







\bibitem{GP09} {\it Perez-Izquierdo J.M.,} Unital algebras, ternary derivations, and local triality, Algebras, Representations and Applications, Contemp. Math., 483, Amer. Math. Soc., Providence, RI, 2009, 205--220

\bibitem{poj03}  {\it Pojidaev A.,} On simple Filippov superalgebras of type $B(0,n)$, J. Algebra Appl., 2 (2003), 3, 335--349.

\bibitem{poj06} {\it Pojidaev A.,  Saraiva P.,} On derivations of the ternary Malcev algebra $M_8$, Comm. Algebra, 34 (2006), 10, 3593--3608.

\bibitem{poj09} {\it Pojidaev A.,  Saraiva P.,} On simple Filippov superalgebras of type $B(0,n)$. II, Port. Math., 66 (2009), 1, 115--130.


\bibitem{poj08}  {\it Pozhidaev A.,} Simple Filippov superalgebras of type $B(m,n)$, Algebra and Logic, 47 (2008), 2, 139--152.


\bibitem{lesha12} {\it Shestakov A.},
Ternary derivations of separable associative and Jordan algebras,
Siberian Math. J., 53 (2012), 5, 943--956.

\bibitem{lesha14} {\it Shestakov A.},
Ternary derivations of Jordan superalgebras,
Algebra and Logic, 53 (2014), 4, 323--348.



\bibitem{Williams}
{\it Williams M. P.,}
Nilpotent n-Lie algebras,
Comm. Alg., 37 (2009), 6, 1843--1849.


\bibitem{ZZ10} {\it Zhang R., Zhang Y., }
Generalized derivations of Lie superalgebras, Comm. Algebra, 38:10 (2010), 3737--3751


\bibitem{ZT15} {\it Zhang T.,}
Cohomology and deformations of 3-Lie colour algebras,
Linear Multilinear Algebra 63 (2015), no. 4, 651--671.

\bibitem{zus10} {\it Zusmanovich~P.,}  On $\delta$-derivations of Lie algebras and superalgebras,
J. of Algebra, 324 (2010), 12, 3470--3486.








\end{thebibliography}
\end{document}